# Symmetry Reduction in Optimal Control of Multi-Agent Systems on Lie Groups

Leonardo J. Colombo and Dimos V. Dimarogonas *Senior Member, IEEE*.

*Abstract*—We study the reduction of degrees of freedom for the equations that determine necessary optimality conditions for extrema in an optimal control problem for a multi-agent system by exploiting the physical symmetries of agents, where the kinematics of each agent is given by a left-invariant control system. Reduced optimality conditions are obtained using techniques from variational calculus and Lagrangian mechanics. A Hamiltonian formalism is also studied, where the problem is explored through an application of Pontryagin's maximum principle for left-invariant systems, and the optimality conditions are obtained as integral curves of a reduced Hamiltonian vector field. We apply the results to an energy-minimum control problem for multiple unicycles.

## I. INTRODUCTION

Dimensionality reduction for large scale systems has become an active problem of interest within the automatic control and robotics communities. In large robotic swarms, guidance and trajectory planning algorithms for coordination while optimizing qualitative features for the swarm of multiple robots are determined by solutions of nonlinear equations which demand a high-computational costs along its integration. The construction of methods for reduction of dimensionality permits fast computations for the generation of optimal trajectories.

Methods for trajectory tracking and estimation algorithms for pose and attitude of mechanical systems evolving on Lie groups are commonly employed for improving accuracy on simulations, as well as to avoid singularities by working with coordinate-free expressions in the associated Lie algebra of the Lie group to describe behaviors in multi-agent systems [13] (i.e., a set of equations depending on an arbitrary choice of the basis for the Lie algebra). More recently, this framework has been used for cooperative transportation [19].

Optimization problems on Lie groups have a long history [14] and have been applied to many problems in control engineering. In practice, many robotic systems exhibit symmetries that can be exploited to reduce some of the complexities of system models, for instance degrees of freedom. Symmetries in optimal control for systems on Lie groups have been studied in [3], [18], [17], [24] among many others, mainly for applications in robotic and aerospace engineering, and in particular, for spacecraft attitude control and underwater vehicles [21]. Recent studies include dynamic programming [22], interconnected systems [11] and soft robotics [6]. While most of the applications of symmetry reduction provided in the literature focus on the single agent situation, only a few works studied the relation between multi-agent systems and symmetry reduction (see for instance the early work on the topic [15]), in this work we introduce a new application in optimal control, by extending the reduction technique to multi-agent systems modeled by left-invariant control systems with a decentralized communication topology determined by an undirected graph, i.e., the information between the agents is only shared between nearest neighbors. Such a method is proposed from a Lagrangian and a Hamiltonian point of view. In particular, in this paper we develop a mathematical framework to reduce the degrees of freedom of the differential equations governing the necessary conditions for optimality in an optimal control problem, and hence as a consequence we are able to achieve reduction in computational costs for the integration of the corresponding equations describing the reduced optimality conditions.

Multi-agent control systems with agents modeled by Lagrangian systems have been recently considered in [1], [7], [8], [9] to study the dynamics and the construction of geometric integrators in formation control and also optimization problems in Riemannian manifolds. In this paper agents are modeled by kinematic left-invariant control systems on a Lie group. These models include nonholonomic, underactuated, fully-actuated, homogeneous and heterogeneous agents. Collision avoidance guarantees are given by introducing appropriate potential functions corresponding to fictitious forces (as for instance, a Coulomb or force field potential) into the cost functional for the optimal control problem according to [16], [20], [23]. Reduced optimality (necessary) conditions for extrema are described via left-trivialized Euler-Lagrange and left-trivialized Hamilton equations.

In our previous work [7] we studied symmetry reduction in optimal control problems for the case of centralize communications between the agents. In this work we extend and improve the results of [7] by (i) considering a decentralized communication between the agents, and (ii) improving the derivation of the equations describing necessary conditions for optimality from variational principles in the Lagrangian framework. In particular, by considering a change of indexes in the cost functional in a step of the proof and by employing

L. Colombo is with Instituto de Ciencias Matemáticas (CSIC-UAM-UC3M-UCM), Calle Nicolas Cabrera 13-15, 28049, Madrid, Spain. leo.colombo@icmat.es

D. Dimarogonas is with Department of Automatic Control, EECS School, KTH Royal Institute of Technology, Osquldas vag 10, 10044, Stockholm, Sweden. dimos@kth.se.

This work was supported by the MINECO (Spain) Grant MTM2016-76702-P, I-Link Project (Ref: linkA20079) from CSIC; "la Caixa" Foundation (ID Project 100010434, code LCF/BQ/PI19/11690016), Swedish Research Council (VR), Knut och Alice Wallenberg foundation (KAW), the H2020 Project Co4Robots and the H2020 ERC Starting Grant BUCOPHSYS.

the neighboring relation between agents we are able to obtain a more compact set of reduced equations in comparison with [7] which allows fast computations from a numerical point of view. Moreover, (iii) we reduced the hypothesis concerning the relations between agents and we obtained a more clear result about the reduced necessary conditions for optimality. Further original contributions of this work (iv) is given by considering a Hamiltonian formalism and reduce the Pontryagin maximum principle to get reduced equations of motion together with an analysis of the relation between both frameworks, the Lagrangian and the Hamiltonian. Finally, in comparison with [7] we improved the results obtained in the proposed example by (v) noting that the coupling term between the agents only depends on the distances of the center of masses between them and therefore, the coupling term is invariant under rotations and translations.

The first characterization of the reduced equations derived in this paper is given by considering the optimal control problem as a constrained problem and impose constraints on the unactuated configurations with Lagrange multipliers. Under symmetry conditions in the cost functions we are able to reduce the Lagrangian associated with the problem by the symmetry group. Reduced optimality conditions are obtained via calculus of variations. The results are given in Theorem 4.1 and Proposition 4.2. The Hamiltonian characterization is given by reducing the Hamiltonian for the optimal control problem via the Pontryagin Maximum Principle (PMP) for left-invariant systems, under the same symmetries conditions. Reduced optimality conditions are obtained as integral curves of a Hamiltonian vector field associated to a reduced Hamiltonian function. The result is given in Theorem 5.1.

While the first method may have useful advantages in practice such as the construction of variational integrators for estimations, the second one provides reduced optimality conditions in an intrinsic way, i.e., independent of the choice of coordinates, in terms of geometric elements. Both approaches provide coordinate-free equations on vector spaces. This is useful in practice because instead of solving differential equations on manifolds which usually requires embedding the manifold into an Euclidean space, and thus increasing the dimension of the system with additional constraints (note that differentiable manifolds are modeled locally on normed spaces and Lie groups are basic examples of differential manifolds), in this work we provide reduced equations on vector spaces. This allows the use of fast and accurate popular numerical solvers for the reduced equations. This is possible because the position and velocity of each agent can always be described as an element in the tangent space at the identity element of the group, by left-translating the velocity vectors based at any point in the group. The tangent space at the identity element of the group corresponds with the Lie algebra of the group, which has a structure of vector space.

It is well known that (see for instance [14]) Hamilton equations (in the cotangent bundle), are the dual representation of Euler-Lagrange equations (in the tangent bundle). Our motivation to develop a Lagrangian formalism for the reduction in optimal control problems is mainly based on the fact that by considering a Lagrangian formalism it is possible to construct variational integrators. That is, a class of geometric numerical schemes that preserves the qualitative features of the system such as momentum preservation and symplecticity, and exhibit a good behavior of the energy preserved along the motion in conservative systems. This can be achieve by discretizing the variational principle, instead of discretizing the equations of motion as is usual in the literature to construct numerical methods for this class of problems. Moreover, it is also well known that Noether's theorem (given in the Lagrangian framework) provides a direct link between symmetries and conserved quantities which is preserved by the discretization of variational principles in the Lagrangian framework.

The structure of the paper is as follows: Section II reviews mechanics on Lie groups. Section III introduces the left-invariant kinematic multi-agent control system and the formulation of the optimal control problem. Sections IV and V are devoted to deriving reduced optimality conditions from a Lagrangian and Hamiltonian point of view, respectively. An example to the optimal control problem of three unicycles is studied in Section VI.

## II. Preliminaries

In this section we introduce the review material we will use along the work. For a further covering of the topics see [3] (Chapter 3) and [12] (Chapters 6-7).

Let $Q$ be the configuration space of a mechanical system, a differentiable manifold of dimension $n$ with local coordinates $q = (q^1, \ldots, q^n)$. Let $TQ$ be the tangent bundle of $Q$, locally described by positions and velocities for the system $v_q = (q^1, \ldots, q^n, \dot{q}^1, \ldots, \dot{q}^n) \in TQ$ with $\dim(TQ) = 2n$. Let $T^*Q$ be its cotangent bundle, locally described by the positions and the momentum for the system, i.e., $(q, p) \in T^*Q$ with $\dim(T^*Q) = 2n$. The tangent bundle at a point $q \in Q$ is denoted as $T_qQ$ and the cotangent bundle at a point $h \in Q$ is denoted as $T_h^*Q$. If the configuration space is a Lie group and the system has a symmetry, one can take advantage of it to reduce the degrees of freedom of the system.

*Definition 2.1:* Let $G$ be a Lie group with identity element $\bar{e} \in G$. A *left-action* of $G$ on a manifold $Q$ is a smooth mapping $\Phi : G \times Q \to Q$ such that $\Phi(\bar{e}, q) = q$ for all $q \in Q$, $\Phi(g, \Phi(h, q)) = \Phi(gh, q)$ for all $g, h \in G, q \in Q$ and for every $g \in G$, $\Phi_g : Q \to Q$ defined by $\Phi_g(q) := \Phi(g, q)$ is a diffeomorphism.

We often use the notation $\Phi_g(q) = \Phi(g, q) := gq$ and say that $g$ acts on $q$. All Lie group actions are assumed to be smooth.

Let $G$ be a finite dimensional Lie group and $\mathfrak{g}$ denotes the Lie algebra associated to $G$ defined as $\mathfrak{g} := T_{\bar{e}}G$, i.e., the tangent space at the identity $\bar{e} \in G$. Let $L_g : G \to G$ be the left translation of the element $g \in G$ given by $L_g(h) = gh$ for $h \in G$. $L_g$ is a diffeomorphism on $G$ and a left-action from $G$ to $G$ [12]. Its tangent maps (i.e, the linearization or tangent lift) is denoted by $T_hL_g : T_hG \to T_{gh}G$. Similarly, the cotangent map (cotangent lift), is defined as $(T_hL_g)^*$, the dual map of the tangent lift denoted by $T_h^*L_g : T_h^*G \to T_{gh}^*G$, and determined by the relation $\langle (T_hL_g)^*(\alpha_{gh}), Y_h \rangle = \langle \alpha_{gh}, (T_hL_g)Y_h \rangle$, $Y_h \in T_hG$, $\alpha_{gh} \in T_{gh}^*G$. It is well known that the tangent and cotangent lift are actions.

Here, $\langle \cdot, \cdot \rangle : \mathfrak{g}^* \times \mathfrak{g} \to \mathbb{R}$ denotes the so-called *natural pairing* between vectors and co-vectors. It is defined by $\langle y, x \rangle := y \cdot x$ for $y \in \mathfrak{g}^*$ and $x \in \mathfrak{g}$, where $x$ is understood as a row vector and $y$ as a column vector. For a matrix Lie algebra $\langle y, x \rangle = y^T x$ (see [12], Section 2.3). Using this pairing between vectors and co-vectors, for $g, h \in G$, $y \in \mathfrak{g}^*$ and $x \in \mathfrak{g}$, one can write

$$\langle T^*_g L_{g^{-1}}(y), T_e L_g(x) \rangle = \langle y, x \rangle. \tag{1}$$

*Definition 2.2:* Denote by $\mathfrak{X}(G)$ the set of vector fields on $G$. A vector field $X \in \mathfrak{X}(G)$ is called *left invariant* if $T_h L_g(X(h)) = X(L_g(h)) = X(gh) \; \forall g, h \in G$.

In particular for $h = \bar{e}$, Definition 2.2 means that a vector field $X$ is left-invariant if $\dot{g} = X(g) = T_{\bar{e}} L_g \xi$ for $\xi = X(\bar{e}) \in \mathfrak{g}$. As $X$ is left invariant, $\xi = X(\bar{e}) = T_g L_{g^{-1}} \dot{g}$. The tangent map $T_{\bar{e}} L_g$ shifts vectors based at $\bar{e}$ to vectors based at $g \in G$. By doing this operation for every $g \in G$ we define a vector field as $\overleftarrow{\xi}(g) := T_{\bar{e}} L_g(\xi)$ for $\xi := X(e) \in T_{\bar{e}} G$. Note that the vector field $\overleftarrow{\xi}(g)$ is left invariant, because $\overleftarrow{\xi}(hg) = T_{\bar{e}}(L_h \circ L_g)\xi = (T_{\bar{e}} L_h) \circ (T_{\bar{e}} L_g)\xi = T_g L_h \overleftarrow{\xi}(g)$. Here we use the upper left-arrow to denote that the vector field is *left* invariant.

*Definition 2.3:* Let $\Phi_g : Q \to Q$ be a left action of $G$ on $Q$; a function $f : Q \to \mathbb{R}$ is said to be *invariant* under the action $\Phi_g$, if $f \circ \Phi_g = f$ (i.e., $\Phi_g$ is a symmetry of $f$).

Consider the vector bundles isomorphisms $\lambda_{TG} : G \times \mathfrak{g} \to TG$ and $\lambda_{T^*G} : G \times \mathfrak{g}^* \to T^*G$ defined as $\lambda_{TG}(g,\xi) = (g, T_{\bar{e}} L_g(\xi))$, $\lambda_{T^*G}(g,\mu) = (g, T_g^* L_{g^{-1}}(\mu))$. $\lambda_{TG}$ and $\lambda_{T^*G}$ are called left-trivializations of $TG$ and $T^*G$ respectively. Therefore, the left-trivialization $\lambda_{TG}$ permits to identify the tangent bundle $TG$ with $G \times \mathfrak{g}$, and through $\lambda_{T^*G}$, the cotangent bundle $T^*G$ can be identified with $G \times \mathfrak{g}^*$.

Denote by $\mathrm{ad} : \mathfrak{g} \times \mathfrak{g} \to \mathfrak{g}$ the *adjoint operator* given by $\mathrm{ad}_\xi \eta := [\xi, \eta]$, where $[\cdot, \cdot]$ denotes the Lie bracket of vector fields on the Lie algebra $\mathfrak{g}$. Let $\mathbf{L} : TG \to \mathbb{R}$ be a Lagrangian function describing the dynamics of a mechanical system. After a left-trivialization of $TG$ we may consider the trivialized Lagrangian $\ell : G \times \mathfrak{g} \to \mathbb{R}$ given by $\ell(g, \xi) = \mathbf{L}(g, T_{\bar{e}} L_{g^{-1}}(\xi)) = \mathbf{L}(g, g\xi)$. The left-trivialized Euler–Lagrange equations on $G \times \mathfrak{g}$ (see, e.g., [12], Ch. 7), are given by

$$\frac{d}{dt} \frac{\partial \ell}{\partial \xi} + T_{\bar{e}}^* L_g \left( \frac{\partial \ell}{\partial g} \right) = \mathrm{ad}^*_\xi \frac{\partial \ell}{\partial \xi}, \tag{2}$$

together with the kinematic equation $\dot{g} = T_{\bar{e}} L_g \xi$. Here, $\mathrm{ad}^* : \mathfrak{g} \times \mathfrak{g}^* \to \mathfrak{g}^*$, $(\xi, \mu) \mapsto \mathrm{ad}^*_\xi \mu$ denotes the *co-adjoint operator*, defined by $\langle \mathrm{ad}^*_\xi \mu, \eta \rangle = \langle \mu, \mathrm{ad}_\xi \eta \rangle$ for all $\eta \in \mathfrak{g}$.

The left-trivialized Euler-Lagrange equations together with the equation $\xi = T_g L_{g^{-1}}(\dot{g})$ are equivalent to the Euler–Lagrange equations for $\mathbf{L}$. Note that for a matrix Lie group, the previous equations is equivalent to $\xi = g^{-1} \dot{g}$.

One can also obtain the reduced Hamiltonian $h : G \times \mathfrak{g}^* \to \mathbb{R}$ given by $h(g, \mu) = \langle \mu, \xi(\mu) \rangle - \ell(g, \xi(\mu))$, where one uses the Legendre transformation $f\ell : G \times \mathfrak{g} \to G \times \mathfrak{g}^*$. If $\ell$ is a diffeomorphism (i.e., $\ell$ is hyper-regular, see [12], Section 4.2, pp. 142), we can define the velocity $\xi$ as a function of the momentum $\mu$ by the implicit function theorem (see [12], Section 9.1, pp. 296 for details in the procedure). The left-trivialized Euler-Lagrange equations (2) can then be written as the left-trivialized Hamilton equations (see, [3], [12]), which are given by $\dot{\mu} = \mathrm{ad}^*_{\frac{\partial h}{\partial \mu}} \mu - T_{\bar{e}}^* L_g \left( \frac{\partial \ell}{\partial g} \right)$.

## III. LEFT-INVARIANT KINEMATIC MULTI-AGENT CONTROL SYSTEM AND PROBLEM FORMULATION

Left-invariant control systems provide a general framework for modeling some classes of systems that include, for instance, the mathematical control design for spacecraft and underwater vehicles (see [21] and references therein). In general, the configuration space for these systems is globally described by a matrix Lie group making it a natural model for a controlled system. This framework gives rise to coordinate-free expressions for the dynamics describing the behavior of the system (i.e., only depends on an arbitrary choice of the basis for the Lie algebra).

### A. Left-invariant kinematic multi-agent control system

Consider a set $\mathcal{N}$ consisting of $r$ free agents evolving each one on a Lie group $G$ with dimension $n$. Along this work we assume that the configuration space of each agent has the same Lie group structure. We note however that each agent can have different masses and inertia values, and therefore agents are heterogeneous. We denote by $g_i \in G$ the configuration (positions) of an agent $i \in \mathcal{N}$ and $g_i(t) \in G$ describes the evolution of agent $i$ at time $t$. The element $g \in G^r$ denotes the stacked vector of positions where $G^r := \underbrace{G \times \ldots \times G}_{r-times}$ denotes the cartesian product of $r$ copies of $G$. We also consider $\mathfrak{g}^r := T_{\bar{e}} G^r$ the Lie algebra associated with the Lie group $G^r$ where $\bar{e} = (\bar{e}_1, \ldots, \bar{e}_r) \in G^r$ is the identity element and $\bar{e}_j$ the identity element of the $j^{th}$-Lie group which determines $G^r$.

When it is necessary we will write $G_i$ and $\mathfrak{g}_i$ to denote the $i^{th}$-Lie group and $i^{th}$-Lie algebra which determines $G^r$ and $\mathfrak{g}^r$, respectively.

The neighbor relationships are described by an undirected graph $\mathcal{G} = (\mathcal{N}, \mathcal{E})$, static and connected, where the set $\mathcal{N}$ describes the vertices of the graph, and where each vertex $i \in \mathcal{N}$ is a left invariant control system, that is, the kinematics of each agent is determined by

$$\dot{g}_i = T_{\bar{e}_i} L_{g_i}(u_i), \quad g_i(0) = g_0^i, \tag{3}$$

where $g_i(\cdot) \in C^1([0, T], G_i)$, $T \in \mathbb{R}$ fixed, and $u_i$, the control input, $u_i = [u_i^1 \ldots u_i^m]^T$ with $m \leq n$ is a curve on the Lie algebra $\mathfrak{g}_i$ of $G_i$. Note that $m \leq n$ so that the control systems we consider can be underactuated or fully actuated.

Given that for all $i \in \mathcal{N}$ the Lie group is the same, we consider that for all the agents $\mathfrak{g} = \mathrm{span}\{e_1, \ldots, e_n\}$, and then $u_i$ is given by $u_i(t) = e_0 + \sum_{k=1}^{m} u_i^k(t) e_k$, where $e_0 \in \mathfrak{g}$. Therefore (3) gives rise to the *kinematic left invariant control systems*

$$\dot{g}_i(t) = g_i(t) \left( e_0 + \sum_{k=1}^{m} u_i^k(t) e_k \right). \tag{4}$$

The set $\mathcal{E} \subset \mathcal{N} \times \mathcal{N}$ denotes the set of edges for $\mathcal{G}$. The set of neighbors for agent $i$ is defined by $\mathcal{N}_i = \{j \in \mathcal{N} : (i, j) \in \mathcal{E}\}$.

Given that our application in Section VI is based on a drift-free system, from now on in the paper, we will consider that the drift term $e_0$ is zero.

*Remark 3.1:* By considering left-invariant systems, agents can exhibit a different number of control inputs, nevertheless, given that we consider the same basis for the Lie algebra we will need to assume an equal number of actuators. Considering a different number of actuators in this framework will need the specification of a supra index $i$ in the element of the basis, that is, $e_k^i$ since each basis will not be the same. Different actuators are a straightforward extension for the results of this work. We decided omit such extension to keep the notation and exposition of the paper at reader-friendly levels. ◇

### B. Problem formulation

Denote by $\pi_i : G^r \to G_i$, $\tau_i : \mathfrak{g}^r \to \mathfrak{g}_i$, $\chi_i : (\mathfrak{g}^*)^r \to \mathfrak{g}_i^*$, $\alpha_i : (T^*G)^r \to T^*G_i$ and $\beta_i : TG^r \to TG_i$ the canonical projections from $G^r$, $\mathfrak{g}^r$, $(\mathfrak{g}^*)^r$, $(T^*G)^r$ and $TG^r$, respectively, over its $i^{th}$-factor.

We want to find optimality conditions in an optimal control problem for the left-invariant multi-agent control system (4), where along their trajectory not only minimize the cost function for the complete networked system, but also ensure that agents avoid collisions with each other.

Similarly as in [20], [23] we assume that each agent $i$ occupies a disk of radius $\bar{r}$ on $G$. The quantity $\bar{r}$ is chosen to be small enough so that it is possible to pack $r$ disks of radius $\bar{r}$ on $G$. We say that agents $i$ and $j$ avoid mutual collision if $||\pi_i(g(t)) - \pi_j(g(t))||_G > \bar{r}$ for all $t$, where $||\cdot||_G$ denotes the norm on $G$. For instance, for matrix Lie groups we can use the Frobenius norm.

If agents $i$ and $j$ avoids mutual collision at initial states, the collision avoidance task in the optimal control problem is guaranteed by introducing distributed collision avoidance potential functions corresponding to fictitious forces into the cost functional. We introduce the potential function (in the sprit of the artificial potential functions [20] and the structural potentials in [23]) $V_{ij} : G \times G \to \mathbb{R}$ with $i \in \mathcal{N}$; $j \in \mathcal{N}_i$, satisfying $V_{ij} = V_{ji}$, and assume these are sufficiently regular.

The problem studied in this work consists on finding reduced optimality conditions in an optimal control problem, taking advantage of the symmetries in the cost functional.

**Problem:** Find optimality (necessary) conditions on the configurations $g(t) = (g_1(t), \ldots, g_r(t))$ and control inputs $u(t) = (u_1(t), \ldots, u_r(t))$ minimizing the cost functional

$$\min_{(g(\cdot), u(\cdot))} \sum_{i=1}^r \int_0^T [C_i(g_i(t), u_i(t)) + \frac{1}{2} \sum_{j \in \mathcal{N}_i} V_{ij}(g_i(t), g_j(t))] dt \quad (5)$$

subject to $\dot{g}_i(t) = T_{\bar{e}_i} L_{g_i(t)}(u_i(t))$ and boundary values $g(0) = (g_1(0), \ldots, g_r(0)) =: (g_1^0, \ldots, g_r^0)$, $g(T) = (g_1(T), \ldots, g_r(T)) =: (g_1^T, \ldots, g_r^T)$, with $u(t) = (u_1(t), \ldots, u_r(t)) \in \mathfrak{g}^r$, and where the cost functions are invariant under the left-action $\rho^i : G \times (G \times \mathfrak{g}) \to G \times \mathfrak{g}$, $\rho_g^i(g_i, u_i) = (L_g g_i, u_i)$, $g \in G$, that is, $C_i \circ \rho_g^i = C_i$.

*Remark 3.2:* Note that the factor $\frac{1}{2}$ in the potential function in (5) comes from the fact that $V_{ij} = V_{ji}$. The cost functions $C_i$ are not related to collision avoidance between agents but only to the energy minimization of each agent. The potential functions used to avoid collision in the proposed approach are essentially distributed collision avoidance potentials as for instance of the kind employed in [16], [20], [23], [10]. ◇

## IV. REDUCED OPTIMALITY CONDITIONS

As in [4], and [17], the optimal control problem can be solved as a constrained variational problem by considering the Lagrange multipliers $\lambda_{g_i} = T_{g_i}^* L_{g_i^{-1}}(\lambda_i(t)) \in T_{g_i}^* G$ with $\lambda_i \in C^1([0, T], \mathfrak{g}^*)$ into the cost functional. The existence of $\lambda_i$ is guaranteed by the Lagrange multiplier Theorem [3]. Let $\mathfrak{g}^* = \text{span}\{e^1, \ldots, e^m, e^{m+1}, \ldots, e^n\}$, where $\{e^1, \ldots, e^n\}$ is the dual basis of the basis $\{e_1, \ldots, e_n\}$ for $\mathfrak{g}$, then $\lambda_i = \sum_{k=m+1}^{n} \lambda_i^k e^k$.

Define the function $\mathbf{C} : G^r \times \mathfrak{g}^r \to \mathbb{R}$, $\mathbf{C}(g, u) = \sum_{i=1}^r C_i(\pi_i(g(t)), \tau_i(u(t)))$, and the extended Lagrangian $\mathcal{L} : G^r \times \mathfrak{g}^r \times (T^*G)^r \to \mathbb{R}$ by

$$\mathcal{L}(g, u, \lambda_g) = \mathbf{C}(g, u) + \sum_{i=1}^r \Big( \langle \alpha_i(\lambda_g(t)), \beta_i(T_{\bar{e}} L_g u(t)) \rangle + \frac{1}{2} \sum_{j \in \mathcal{N}_i} V_{ij}(\pi_i(g(t)), \pi_j(g(t))) \Big).$$

Consider the left-action $\rho : G \times (G^r \times \mathfrak{g}^r) \to G^r \times \mathfrak{g}^r$, $\rho_h(g, u) = (L_h g_1, \ldots, L_h g_r, u)$, $h \in G$. The following result gives rise to reduced optimality (necessary) conditions for extremals in the optimal control problem.

*Theorem 4.1:* If $\mathbf{C}$ is invariant under the left action $\rho_g : G^r \times \mathfrak{g}^r \to G^r \times \mathfrak{g}^r$, extremals of the cost functional for the problem (5) satisfy the equations

$$0 = \frac{d}{dt}\left(\frac{\partial C_i}{\partial u_i} + \lambda_i\right) - \text{ad}_{u_i}^*\left(\frac{\partial C_i}{\partial u_i} + \lambda_i\right) - \sum_{j \in \mathcal{N}_i} T_{\bar{e}_i}^* L_{g_i}\left(\frac{\partial \mathbf{V}_{ij}}{\partial g_i}\right), \quad (6)$$

together with $\dot{g}_i = T_{\bar{e}_i} L_{g_i}(u_i)$ for $i = 1, \ldots, r$, where $\mathbf{V}_{ij} : G \to \mathbb{R}$ is given by $V_{ij}(\bar{e}_i, g_i^{-1} g_j)$.

*Proof of Theorem 4.1:* Given that $\lambda_{g_i} = T_{g_i}^* L_{g_i^{-1}}(\lambda_i)$ and $(T_{g_i} L_{g_i^{-1}} \circ T_{\bar{e}_i} L_{g_i}) = \bar{e}_i$, by using (1) we have,

$$\langle \lambda_{g_i}, T_{\bar{e}_i} L_{g_i} u_i \rangle = \langle T_{g_i}^* L_{g_i^{-1}}(\lambda_i), T_{\bar{e}_i} L_{g_i} u_i \rangle = \langle \lambda_i, u_i \rangle.$$

The invariance of $\mathbf{C}$ under $\rho_g$, makes possible to define the reduced Lagrangian $\ell : G^{r-1} \times \mathfrak{g}^r \times (\mathfrak{g}^*)^r \to \mathbb{R}$ as

$$\ell(g, u, \lambda) = \sum_{i=1}^r (C_i(\tau_i(u(t))) + \langle \chi_i(\lambda), \tau_i(u) \rangle + \frac{1}{2} \sum_{j \in \mathcal{N}_i} \mathbf{V}_{ij}(\pi_i(g^{-1}(t)) \pi_j(g(t)))).$$

where $\mathbf{V}_{ij} : G \to \mathbb{R}$ is given by $V_{ij}(\bar{e}_i, g_i^{-1} g_j)$. Note that here $\ell : G^{r-1} \times \mathfrak{g}^r \times (\mathfrak{g}^*)^r \to \mathbb{R}$, nevertheless with a slight abuse of the notation we denote by $g \in G^{r-1}$ one of the inputs for $\ell$, whereas previously $g$ was an element of $G^r$. In what follows $g$ should be considered as an element of $G^{r-1}$ instead of $G^r$ in $\ell$ and an element of $G^r$ in $\mathcal{L}$.

Next, after obtaining the reduced Lagrangian we shown that for variations of $g$ vanishing at end points, that is, $\delta g(0) = \delta g(T) = 0$, and variations of $u$, the variational principle

$$\delta \int_0^T \mathcal{L}(g(t), u(t), \lambda_g(t)) dt = 0 \qquad (7)$$

implies the constrained variational principle

$$\delta \int_0^T \ell(g(t), u(t), \lambda(t)) dt = 0 \qquad (8)$$

for variations $\delta u = \dot{\eta} + \mathrm{ad}_u \eta$ where $\eta(t) \in C^1([0,T], \mathfrak{g}^r)$, $\eta(0) = \eta(T) = 0$.

Using $\mathbf{C} \circ \rho_g = \mathbf{C}$ and $\langle \lambda_{g_i}, T_{\bar{e}_i} L_{g_i} u_i \rangle = \langle \lambda_i, u_i \rangle$, both integrands are equal, and the variations of $g_i$, $\delta g_i$, induce and are induced by variations $\delta u_i = \dot{\eta}_i + \mathrm{ad}_{u_i} \eta_i$ with $\eta_i(0) = \eta_i(T) = 0$ (See [12] Section 7.3, pp 255). Therefore, if we choose variations such that $\eta_i = T_{g_i} L_{(g_i)^{-1}}(\delta g_i)$ (that is, $\delta g_i = g_i \eta_i$), $\delta u_i = T_{g_i} L_{(g_i)^{-1}}(\dot{g}_i)$ and the variational principle (7) holds, it follows that $\delta u_i = \dot{\eta}_i + \mathrm{ad}_{u_i} \eta_i$ and hence the variational principle (7) implies the constrained variational principle (8).

Now, note that

$$\delta \int_0^T \ell(g(t), u(t), \lambda(t)) dt = \sum_{i=1}^r \int_0^T \left( \left\langle \frac{\partial C_i}{\partial u_i}, \delta u_i \right\rangle + \langle \lambda_i, \delta u_i \rangle \right.$$
$$\left. + \frac{1}{2} \sum_{j \in \mathcal{N}_i} \left\langle \frac{\partial \mathbf{V}_{ij}}{\partial g_i}, \delta g_i \right\rangle + \frac{1}{2} \sum_{j \in \mathcal{N}_i} \left\langle \frac{\partial \mathbf{V}_{ij}}{\partial g_j}, \delta g_j \right\rangle \right) dt$$
$$= \sum_{i=1}^r \int_0^T \left( \left\langle \frac{\partial C_i}{\partial u_i} + \lambda_i, \dot{\eta}_i + \mathrm{ad}_{u_i} \eta_i \right\rangle \right.$$
$$\left. + \frac{1}{2} \sum_{j \in \mathcal{N}_i} \left\langle \frac{\partial \mathbf{V}_{ij}}{\partial g_i}, \delta g_i \right\rangle + \frac{1}{2} \sum_{j \in \mathcal{N}_i} \left\langle \frac{\partial \mathbf{V}_{ij}}{\partial g_j}, \delta g_j \right\rangle \right) dt.$$

where the first equality comes from the definition of variation of a function on a manifold, that is, $\delta f(\xi) = \langle \frac{\partial f}{\partial \xi}, \delta \xi \rangle$ for an arbitrary function $f$ in an arbitrary manifold, and the second one by replacing the variations by their corresponding expressions given before.

The first component of the previous integrand, after applying integration by parts twice, using the boundary conditions and the definition of co-adjoint action, results in

$$\sum_{i=1}^r \int_0^T \left\langle -\frac{d}{dt}\left(\frac{\partial C_i}{\partial u_i} + \lambda_i\right) + \mathrm{ad}^*_{u_i}\left(\frac{\partial C_i}{\partial u_i} + \lambda_i\right), \eta_i \right\rangle dt.$$

Using the fact that $T(L_{g_i} \circ L_{g_i^{-1}}) = TL_{g_i} \circ TL_{g_i^{-1}}$ is equal to the identity map on $TG$ and $\eta_i = T_{g_i} L_{g_i^{-1}}(\delta g_i)$, the second component can be written as

$$\sum_{j \in \mathcal{N}_i} \left\langle \frac{\partial \mathbf{V}_{ij}}{\partial g_i}, \delta g_i \right\rangle = \sum_{j \in \mathcal{N}_i} \left\langle \frac{\partial \mathbf{V}_{ij}}{\partial g_i}, (T_{\bar{e}_i} L_{g_i} \circ T_{g_i} L_{g_i^{-1}}) \delta g_i \right\rangle$$
$$= \sum_{j \in \mathcal{N}_i} \left\langle \frac{\partial \mathbf{V}_{ij}}{\partial g_i}, T_{\bar{e}_i} L_{g_i} \cdot \eta_i \right\rangle$$
$$= \sum_{j \in \mathcal{N}_i} \left\langle T^*_{\bar{e}_i} L_{g_i} \left(\frac{\partial \mathbf{V}_{ij}}{\partial g_i}\right), \eta_i \right\rangle$$

For the last member of the integrand, using the fact that $V_{ij} = V_{ji}$, we get $\mathbf{V}_{ij}(g_i^{-1} g_j) = V_{ij}(\bar{e}_i, g_i^{-1} g_j) = V_{ji}(\bar{e}_i, g_i^{-1} g_j) = \mathbf{V}_{ji}(g_i^{-1} g_j)$, and by employing a change of variables, it can be written as $\sum_{j \in \mathcal{N}_i} \left\langle T^*_{\bar{e}_i} L_{g_i} \frac{\partial \mathbf{V}_{ij}}{\partial g_i}, \eta_i \right\rangle$.

Therefore, $\delta \int_0^T \ell(g(t), u(t), \lambda(t)) dt = 0$, $\forall \delta \eta_i$ implies

$$0 = \frac{d}{dt}\left(\frac{\partial C_i}{\partial u_i} + \lambda_i\right) - \mathrm{ad}^*_{u_i}\left(\frac{\partial C_i}{\partial u_i} + \lambda_i\right) - \sum_{j \in \mathcal{N}_i} T^*_{\bar{e}_i} L_{g_i}\left(\frac{\partial \mathbf{V}_{ij}}{\partial g_i}\right).$$

Note that the previous equations are on Lie algebras but with a coupled term on $G$. To describe the dynamics of each individual agent into the Lie group and therefore obtain the configurations $g(t)$, together with the Lie algebra elements $u_i \in \mathfrak{g}$ and $\lambda_i \in \mathfrak{g}^*$, we need also to consider the equation $\dot{g}_i = T_{\bar{e}_i} L_{g_i}(u_i)$ and solve the coupled system of equations. □

Equations (6) can not describe completely the time evolution of the controls and the Lagrange multipliers. Since they are two independent variables, we must have two equations in order to have a system of differential equations with a well defined solution. To tackle this issue we propose the following splitting of the equations.

*Proposition 4.2:* If the Lie algebra admits a decomposition of the form $\mathfrak{g} = \mathfrak{r} \oplus \mathfrak{s}$, where $\mathfrak{r} = \mathrm{span}\{e_1, \ldots, e_m\}$, $\mathfrak{s} = \mathrm{span}\{e_{m+1}, \ldots, e_n\}$, such that

$$[\mathfrak{s}, \mathfrak{s}] \subseteq \mathfrak{s}, \quad [\mathfrak{s}, \mathfrak{r}] \subseteq \mathfrak{r}, \quad [\mathfrak{r}, \mathfrak{r}] \subseteq \mathfrak{s}, \qquad (9)$$

then the time evolution of equations (6) can be rewritten as

$$\frac{d}{dt} \frac{\partial C_i}{\partial u_i} = \mathrm{ad}^*_{u_i} \lambda_i \Big|_{\mathfrak{r}} + \sum_{j \in \mathcal{N}_i} T^*_{\bar{e}_i} L_{g_i}\left(\frac{\partial \mathbf{V}_{ij}}{\partial g_i}\right)\Big|_{\mathfrak{r}}, \qquad (10)$$

$$\dot{\lambda}_i = \mathrm{ad}^*_{u_i} \frac{\partial C_i}{\partial u_i} \Big|_{\mathfrak{s}} + \sum_{j \in \mathcal{N}_i} T^*_{\bar{e}_i} L_{g_i}\left(\frac{\partial \mathbf{V}_{ij}}{\partial g_i}\right)\Big|_{\mathfrak{s}}. \qquad (11)$$

*Proof:* Given that $\mathfrak{g} = \mathfrak{r} \oplus \mathfrak{s}$ it follows that $\mathfrak{g}^* = \mathfrak{r}^* \oplus \mathfrak{s}^*$ where $\mathfrak{r}^* = \mathrm{span}\{e^1, \ldots, e^m\}$ and $\mathfrak{s}^* = \mathrm{span}\{e^{m+1}, \ldots, e^n\}$. Moreover, using (9), this last decomposition satisfies that $\mathrm{ad}^*_{\mathfrak{r}} \mathfrak{s}^* \subseteq \mathfrak{r}^*$, $\mathrm{ad}^*_{\mathfrak{s}} \mathfrak{r}^* \subseteq \mathfrak{r}^*$, $\mathrm{ad}^*_{\mathfrak{s}} \mathfrak{s}^* \subseteq \mathfrak{s}^*$, and $\mathrm{ad}^*_{\mathfrak{r}} \mathfrak{r}^* \subseteq \mathfrak{s}^*$, and therefore, given that, $\frac{\partial C_i}{\partial u_i} \in \mathfrak{r}^*$ and $\lambda_i \in \mathfrak{s}^*$ we have $\mathrm{ad}^*_{u_i}\left(\frac{\partial C_i}{\partial u_i}\right) \in \mathfrak{s}^*$ and $\mathrm{ad}^*_{u_i} \lambda_i \in \mathfrak{r}^*$. Using the previous decomposition, the second factor in the right hand side of (6) can be split into $\mathrm{ad}^*_{u_i}\left(\frac{\partial C_i}{\partial u_i}\right)$ and $\mathrm{ad}^*_{u_i} \lambda_i$.

Since $T^*_{\bar{e}_i} L_{g_i}\left(\frac{\partial \mathbf{V}_{ij}}{\partial g_i}\right) \in \mathfrak{g}^*$, it has a decomposition into $\mathfrak{r}^*$ and $\mathfrak{s}^*$ as

$$\sum_{k=1}^m \left(T^*_{\bar{e}_i} L_{g_i}\left(\frac{\partial \mathbf{V}_{ij}}{\partial g_i}\right)\right) e_k \in \mathfrak{r}^*, \quad \sum_{k=m+1}^n \left(T^*_{\bar{e}_i} L_{g_i}\left(\frac{\partial \mathbf{V}_{ij}}{\partial g_i}\right)\right) e_k \in \mathfrak{s}^*.$$

Hence, for all $i$ (6) splits into the following equations

$$\frac{d}{dt} \frac{\partial C_i}{\partial u_i} = \mathrm{ad}^*_{u_i} \lambda_i \Big|_{\mathfrak{r}} + \sum_{j \in \mathcal{N}_i} T^*_{\bar{e}_i} L_{g_i}\left(\frac{\partial \mathbf{V}_{ij}}{\partial g_i}\right)\Big|_{\mathfrak{r}},$$

$$\dot{\lambda}_i = \mathrm{ad}^*_{u_i} \frac{\partial C_i}{\partial u_i}\Big|_{\mathfrak{s}} + \sum_{j \in \mathcal{N}_i} T^*_{\bar{e}_i} L_{g_i}\left(\frac{\partial \mathbf{V}_{ij}}{\partial g_i}\right)\Big|_{\mathfrak{s}}.$$

□

*Remark 4.3:* Note that if the reduced Lagrangian $\ell$ is hyper-regular (i.e., $\ell$ is a diffeomorphism - see Section 2), then given initial conditions for the system (3)-(10)-(11), solutions of the initial value problem exist and are unique. ◊

*Remark 4.4:* Note that semi-simple Lie algebras admit a Cartan decomposition, i.e., if $\mathfrak{g}$ is semi-simple, then $\mathfrak{g} = \mathfrak{r} \oplus \mathfrak{s}$ such that $[\mathfrak{s}, \mathfrak{s}] \subseteq \mathfrak{s}$, $[\mathfrak{s}, \mathfrak{r}] \subseteq \mathfrak{r}$, $[\mathfrak{r}, \mathfrak{r}] \subseteq \mathfrak{s}$. The converse, however, is not necessarily true. In particular, we do not restrict our analysis to semi-simple Lie algebras. The result given in Proposition 4.2 says that if the Lie algebra $\mathfrak{g}$ admits such a decomposition which implies $\mathfrak{g} = \mathfrak{r} \oplus \mathfrak{s}$ with $[\mathfrak{s}, \mathfrak{s}] \subseteq \mathfrak{s}$, $[\mathfrak{s}, \mathfrak{r}] \subseteq \mathfrak{r}$, $[\mathfrak{r}, \mathfrak{r}] \subseteq \mathfrak{s}$, then we can split the equations, without any further assumption on the structure of the Lie algebra. Moreover, in Section VI, we study an example for a non-semi-simple Lie algebra which yet exhibits such a decomposition. ◊

## V. REDUCED NECESSARY CONDITION VIA THE REDUCED PONTRYAGIN MAXIMUM PRINCIPLE

Next, we show how Hamilton's principle defines an optimal control problem for which an appropriate Hamiltonian $h : G^{r-1} \times (\mathfrak{g}^*)^r \to \mathbb{R}$ is obtained through an application of Pontryagin's maximum principle. Reduced optimality (necessary) conditions for extrema are obtained as integral curves of the Hamiltonian vector field for $h$.

### A. Reduced optimality of conditions via the reduced PMP

Consider the optimal control problem given in (5). A Hamiltonian structure comes into play through the augmented cost functional $\mathcal{S}^a$ defined on the space of smooth functions from $[0, T]$ to $(T^*G)^r \times_{G^r} (G^r \times \mathfrak{g}^r)$ given by

$$\mathcal{S}^a(g, u, \mu_g) = \int_0^T \left( \mathbf{C}(g, u) + \sum_{i=1}^r \langle \alpha_i(\mu_g), \beta_i(\dot{g}) - \beta_i(T_{\bar{e}}L_g u) \rangle \right.$$
$$\left. + \sum_{i=1}^r \sum_{j \in \mathcal{N}_i} V_{ij}(\pi_i(g), \pi_j(g)) \right) dt,$$

where $\alpha_i(\mu_g) = \mu_{g_i}(t) = T^*_{g_i} L_{g_i^{-1}}(\mu_i(t)) \in T^*_{g_i} G$ with $\mu_i \in C^1([0, T], \mathfrak{g}^*)$. We used the notation $\times_{G^r}$ to denote the product bundle with fibers on the manifold $G^r$, this means that in the space $(T^*G)^r \times_{G^r} (G^r \times \mathfrak{g}^r)$, the $g$ component in the space $(T^*G)^r$ is the same as the $g$ component in $(G^r \times \mathfrak{g}^r)$. The augmented cost function permits to introduce the control Hamiltonian $\mathcal{H}_c : (T^*G)^r \times_{G^r} (G^r \times \mathfrak{g}^r) \to \mathbb{R}$ as

$$\mathcal{H}_c(g, \mu_g, u) = -\mathbf{C}(g, u) + \sum_{i=1}^r \langle \alpha_i(\mu_g(t)), \beta_i(T_{\bar{e}}L_g u(t)) \rangle$$
$$- \sum_{i=1}^r \sum_{j \in \mathcal{N}_i} V_{ij}(\pi_i(g), \pi_j(g)).$$

Along the proof for reduced optimality conditions we will employ Pontryagin's maximum principle for left invariant control systems (see Theorem 2.1 in [15] and [14]).

*Theorem 5.1:* If $\mathbf{C}$ is invariant under the left-action $\rho : G \times (G^r \times \mathfrak{g}^r) \to G^r \times \mathfrak{g}^r$, $\rho_h(g, u) = (L_h g_1, \ldots, L_h g_r, u)$, $h \in G$, $u \in \mathfrak{g}^r$, reduced optimality conditions for extrema are determined by integral curves of the Hamiltonian vector field $X_h$ for the reduced Hamiltonian $h : G^{r-1} \times (\mathfrak{g}^*)^r \to \mathbb{R}$ satisfying Hamilton's equations for $h(g, \mu)$

$$\dot{g}_i = T_{\bar{e}_i} L_{g_i} u_i^\star, \quad \dot{\mu}_i = \mathrm{ad}^*_{u_i^\star} \mu_i - \sum_{j \in \mathcal{N}_i} T^*_{\bar{e}_i} L_{g_i} \left( \frac{\partial \mathbf{V}_{ij}}{\partial g_i} \right),$$

where $u^\star = (u_1^\star, \ldots, u_r^\star)$ denotes the optimal control for $\mathcal{H}_c$ and $\mathbf{V}_{ij} : G \to \mathbb{R}$ is given by $V_{ij}(\bar{e}_i, g_i^{-1} g_j)$.

*Proof:* Consider the controlled Hamiltonian $\mathcal{H}_c$. By Pontryagin's maximum principle, we can define the optimal Hamiltonian $\mathcal{H} : (T^*G)^r \to \mathbb{R}$ by

$$\mathcal{H}(g, \mu_g) := \max_{u(\cdot)} \mathcal{H}_c(g, \mu_g, u) = \mathcal{H}_c(g, \mu_g, u^\star)$$

where $u^\star$ denotes the optimal control, determined by the maximization of the Hamiltonian.

Given that each $\mathbf{C} \circ \rho_g = \mathbf{C}$ and $\mu_{g_i}(t) = T^*_{g_i} L_{g_i^{-1}}(\mu_i(t))$, the left action induces the reduced optimal Hamiltonian $h : G^{r-1} \times (\mathfrak{g}^*)^r \to \mathbb{R}$ given by

$$h(g, \mu) = -\mathbf{C}(u^\star) + \sum_{i=1}^r \left( \langle \chi_i(\mu), \tau_i(u^\star) \rangle - \sum_{j \in \mathcal{N}_i} \mathbf{V}_{ij}(\pi_i(g^{-1}) \pi_j(g)) \right)$$

where $\chi_i(\mu) = \mu_i \in \mathfrak{g}_i$ with $\mu_{g_i} = T^*_{g_i} L_{g_i^{-1}}(\mu_i)$ and $\mathbf{V}_{ij}$ is given by $V_{ij}(\bar{e}_i, g_i^{-1} g_j)$. Note that here, as in the Lagrangian case, we are doing an abuse of notation by denoting as $g \in G^{r-1}$ one of the inputs for $h : G^{r-1} \times (\mathfrak{g}^r)^*$, whereas previously $g$ was an element of $G^r$. In what follows $g$ must be considered an element of $G^{r-1}$ instead of $G^r$ for $h$.

Next, we find Hamilton's equations for the reduced optimal Hamiltonian $h$. That is, (by definition) we must find the Hamiltonian vector field $X_h$ solution for $i_{X_h} \omega_T = dh$ where $\omega_T$ is the left trivializations for the canonical symplectic structure on $G^r \times (\mathfrak{g}^*)^r$ given by (See [2])

$$(\omega_T)_{(g,\mu)}((\xi^1, \nu^1), (\xi^2, \nu^2)) = -\langle \nu^1, \xi^2 \rangle + \langle \nu^2, \xi^1 \rangle + \langle \mu, [\xi^1, \xi^2] \rangle \tag{12}$$

where $(\xi^1, \nu^1), (\xi^2, \nu^2) \in \mathfrak{g}^r \times (\mathfrak{g}^*)^r$ and $\xi^i = T_g L_{g^{-1}} v_i$, $v_i \in T_g G^r$ with $i = 1, 2$.

Consider the Hamiltonian vector field $X_h$ for $h$, that is $X_h(g, \mu) = (\xi^1, \mu^1)$, with $\xi^1 \in \mathfrak{g}^r$ and $\mu^1 \in (\mathfrak{g}^*)^r$. By computing the differential of $h$ we obtain

$$dh_{(g,\mu)}(\xi^2, \nu^2) = \sum_{i=1}^r \langle \nu_i^2, u_i^\star \rangle - \sum_{j \in \mathcal{N}_i} T^*_{\bar{e}_i} L_{g_i} \left( \frac{\partial \mathbf{V}_{ij}}{\partial g_i} \right) \tag{13}$$

for $(\xi^2, \nu^2) \in \mathfrak{g}^r \times (\mathfrak{g}^*)^r$ (observe that $\frac{\partial h}{\partial \mu}(g, \mu) \in \mathfrak{g}^{**} = \mathfrak{g}$ since $G$ is finite dimensional).

Note that $i_{X_h}(\omega_T) = \omega_T((\xi^1, \nu^1), (\xi^2, \nu^2))$. Therefore, equating the expression for $\omega_T$ given in (12) with $dh$ given in (13), and using that $\langle \mu, [\xi^1, \xi^2] \rangle = \langle \mu, \mathrm{ad}_{\xi^1} \xi^2 \rangle = \langle \mathrm{ad}^*_{\xi^1} \mu, \xi^2 \rangle$, we obtain $\xi_i^1 = u_i^\star$ and $\nu_i^1 = \mathrm{ad}^*_{u_i^\star} \mu_i + \sum_{j \in \mathcal{N}_i} T^*_{\bar{e}_i} L_{g_i} \left( \frac{\partial \tilde{\mathbf{V}}_{ij}}{\partial g_i} \right)$.

Taking $\dot{g}_i = g_i u_i^\star$ it follows that integral curves for the Hamiltonian vector field $X_h$ must satisfy $\dot{g}_i = T_{\bar{e}_i} L_{g_i} u_i^\star$, and

$$\dot{\mu}_i = \mathrm{ad}^*_{u_i^\star} \mu_i - \sum_{j \in \mathcal{N}_i} T^*_{\bar{e}_i} L_{g_i} \left( \frac{\partial \mathbf{V}_{ij}}{\partial g_i} \right). \quad \square$$

*Remark 5.2:* Under a regularity assumption, it is possible to show the equivalence between the two formalisms presented in the paper, that is, it is possible to obtain the left-trivialized Hamilton equations associated with the optimal control problem from the variational formalism (and vice-versa) by inducing a Legendre transformation. If the reduced Lagrangian $\ell : G^{r-1} \times \mathfrak{g}^r \times (\mathfrak{g}^*)^r \to \mathbb{R}$ is hyper-regular, then the reduced Hamiltonian $h : G^{r-1} \times \mathfrak{g}^r \times (\mathfrak{g}^*)^r \to \mathbb{R}$ induced by the Legendre transformation (which is not the reduced Hamiltonian obtained from the Pontryagin maximum principle) (see [12], Section 9.1, pp 296 for details for the general derivation of $h$ using $f\ell$), after prescribing the optimal controls $u_i^\star$, is given by: $h(g, \mu, \lambda) = \sum_{i=1}^r \langle \mu_i, u_i^\star \rangle - \ell(g, u^\star, \lambda)$, where $\mu_i = \frac{\partial \ell}{\partial u_i^\star} = (\frac{\partial C_i}{\partial u_i^\star} + \lambda_i)$. The left trivialized Euler–Lagrange equations (6) can now be written as the left trivialized Hamilton equations (see Section II), which are given by $\dot{\mu}_i = \mathrm{ad}^*_{u_i^\star} \mu_i - \sum_{j \in \mathcal{N}_i} T^*_{\bar{e}_i} L_{g_i} \left( \frac{\partial \mathbf{V}_{ij}}{\partial g_i} \right)$ for $i = 1, \ldots, r$, together with $\dot{g}_i = T_{\bar{e}_i} L_{g_i} \left( \frac{\partial h}{\partial \mu_i} \right) = T_{\bar{e}_i} L_{g_i} u_i^\star$. ◇

## VI. Optimal Control of Multiples Unicycles

### A. Unicycle model

A unicycle is a homogeneous disk rolling on a horizontal plane maintaining its vertical position (see, e.g. [3]). The configuration of each unicycle at any given time is determined by the element $g_i \in SE(2) \cong \mathbb{R}^2 \times SO(2)$ given by
$$g_i = \begin{bmatrix} \cos\theta_i & -\sin\theta_i & x_i \\ \sin\theta_i & \cos\theta_i & y_i \\ 0 & 0 & 1 \end{bmatrix}, \quad i = 1, 2, 3$$
where $(x_i, y_i) \in \mathbb{R}^2$ represents the point of contact of each wheel with the ground and $\theta_i \in SO(2)$ represents the angular orientation of each unicycle. We denote $u_i = (u_i^1, u_i^2)$. The control input $u_i^1$ represents a force applied to the center of mass of the unicycle and $u_i^2$ a torque applied about its vertical axis. The kinematic equations for the multi-agent system are
$$\dot{x}_i = u_i^2 \cos\theta_i, \quad \dot{y}_i = u_i^2 \sin\theta_i, \quad \dot{\theta}_i = u_i^1, \quad i = 1, 2, 3. \quad (14)$$

### B. Reduction of necessary conditions (Lagrangian)

Equations (14) on $(SE(2))^3$ gives rise to a left-invariant control system where equations take the form $\dot{g}_i = g_i(e_1 u_i^1 + e_2 u_i^2)$ describing all directions of allowable motion, where the elements of the basis of $\mathfrak{se}(2)$ are
$$e_1 = \begin{bmatrix} 0 & -1 & 0 \\ 1 & 0 & 0 \\ 0 & 0 & 0 \end{bmatrix}, \quad e_2 = \begin{bmatrix} 0 & 0 & 1 \\ 0 & 0 & 0 \\ 0 & 0 & 0 \end{bmatrix}, \quad e_3 = \begin{bmatrix} 0 & 0 & 0 \\ 0 & 0 & 1 \\ 0 & 0 & 0 \end{bmatrix},$$
which satisfy $[e_1, e_2] = e_3$, $[e_2, e_3] = 0_{3\times 3}$, $[e_3, e_1] = e_2$. Using the dual pairing, where $\langle \alpha, \xi \rangle := \mathrm{tr}(\alpha \xi)$, for $\xi \in \mathfrak{se}(2)$ and $\alpha \in \mathfrak{se}(2)^*$, the elements of the basis of $\mathfrak{se}(2)^*$ are
$$e^1 = \begin{bmatrix} 0 & \frac{1}{2} & 0 \\ -\frac{1}{2} & 0 & 0 \\ 0 & 0 & 0 \end{bmatrix}, \quad e^2 = \begin{bmatrix} 0 & 0 & 0 \\ 0 & 0 & 0 \\ 1 & 0 & 0 \end{bmatrix}, \quad e^3 = \begin{bmatrix} 0 & 0 & 0 \\ 0 & 0 & 0 \\ 0 & 1 & 0 \end{bmatrix}.$$

Here, $\mathfrak{r} = \{e_1, e_2\}$, $\mathfrak{s} = \{e_3\}$, $\mathfrak{se}(2) = \mathfrak{r} \oplus \mathfrak{s}$ and fulfill the hypothesis of Proposition 4.2. Also note that $\mathfrak{se}(2)$ is not a semi-simple Lie algebra but it satisfies the Lie bracket relations to decompose the dynamics as in Proposition 4.2.

Consider the potential functions $V_{ij} : SE(2) \times SE(2) \to \mathbb{R}$,
$$V_{ij}(g_i, g_j) = \frac{\sigma_{ij}}{2((x_i - x_j)^2 + (y_i - y_j)^2 - d_{ij}^2)},$$
where $\sigma_{ij} \in \mathbb{R}^+$, $d_{ij} \in \mathbb{R}^+$ are prescribed distances between the agents, and also consider the cost functions $C_i(g_i, u_i) = \frac{1}{2} \langle u_i, u_i \rangle$. Since $V_{ij}$ only depends on the distance of the center of masses, it is invariant under rotations and translations (see [23] for instance), i.e., $V_{ij}(hg_i, hg_j) = V_{ij}(g_i, g_j)$ for $h \in SE(2)$.

Denote by
$$\Gamma_{ij} := T^*_{\bar{e}_i} L_{g_i} \left( \frac{\partial \mathbf{V}_{ij}}{\partial g_i} \right) \bigg|_{\mathfrak{r}_i} = \frac{-\sigma_{ij}(x_j - x_i)}{16((x_j - x_i)^2 + (y_j - y_i)^2 - d_{ij}^2)^2} e^2,$$
$$\widetilde{\Gamma}_{ij} := T^*_{\bar{e}_i} L_{g_i} \left( \frac{\partial \mathbf{V}_{ij}}{\partial g_i} \right) \bigg|_{\mathfrak{s}_i} = \frac{-\sigma_{ij}(y_j - y_i)}{16((x_j - x_i)^2 + (y_j - y_i)^2 - d_{ij}^2)^2} e^3.$$

The equations obtained by employing Proposition 4.2 are
$$\dot{u}_i^1 = -\frac{u_i^2 \lambda_i^3}{2}, \quad \dot{u}_i^2 = u_i^1 \lambda_i^3 - \sum_{j \in \mathcal{N}_i} (\Gamma_{ij})_{13}, \quad \dot{\lambda}_i^3 = -u_i^1 u_i^2 + \sum_{j \in \mathcal{N}_i} (\widetilde{\Gamma}_{ij})_{23},$$
together with $\dot{g}_i = g_i(e_1 u_i^1 + e_2 u_i^2)$ for $i = 1, 2, 3$.

### C. Reduction of necessary conditions (Hamiltonian)

The augmented cost functional is given by
$$\mathcal{S}^a(g, \mu_g, u) = \int_0^T \left( \sum_{i=1}^3 \frac{1}{2} \langle u_i, u_i \rangle + \langle T^*_{g_i} L_{g_i^{-1}}(\mu_i), g_i u_i \rangle + \sum_{j \in \mathcal{N}_i} V_{ij}(g_i, g_j) \right) dt.$$
where $\mu_{g_i} = T^*_{g_i} L_{g_i^{-1}}(\mu_i(t)) \in T^*_{\mu_{g_i}} SE(2)$ and $\mu_i \in \mathfrak{se}(2)^*$. The augmented cost functional $\mathcal{S}^a$ induces the control Hamiltonian
$$\mathcal{H}_c(g, \mu_g, u) = \sum_{i=1}^3 \langle T^*_{g_i} L_{g_i^{-1}}(\mu_i), g_i u_i \rangle - \frac{1}{2} \langle u_i, u_i \rangle - \sum_{j \in \mathcal{N}_i} V_{ij}(g_i, g_j).$$
By applying Pontryagin's maximum principle to the control Hamiltonian $\mathcal{H}_c$ we obtain the optimal Hamiltonian $\mathcal{H} : (T^*SE(2))^3 \to \mathbb{R}$ given by
$$\mathcal{H}(g, \mu_g) = \sum_{i=1}^3 \langle T^*_{g_i} L_{g_i^{-1}}(\mu_i), g_i u_i^\star \rangle - \frac{1}{2} \langle u_i^\star, u_i^\star \rangle - \sum_{j \in \mathcal{N}_i} V_{ij}(g_i, g_j)$$
where $u_i^\star$ is the optimal control.

We use the maximization condition $\frac{\partial \mathcal{H}_c}{\partial u_i^\star} = 0$ to write, in the basis of $\mathfrak{se}(2)^*$, the controls $u_i^\star$ in terms of the momenta $\mu_i$, that is, $(u_1^1)^\star = \frac{1}{2} \mu_1^1$, $(u_1^2)^\star = \mu_1^2$, $(u_2^1)^\star = \frac{1}{2} \mu_2^1$ and $(u_2^2)^\star = \mu_2^2$.

The reduced hamiltonian $h : SE(2)^2 \times (\mathfrak{se}(2)^*)^3 \to \mathbb{R}$ is
$$h(\mu, \chi) = \frac{3}{8}((\mu_1^1)^2 + (\mu_1^2)^2) + \frac{1}{2}((\mu_2^1)^2 + (\mu_2^2)^2) - \sum_{j \in \mathcal{N}_i} \mathbf{V}_{ij}(g_i^{-1} g_j).$$

Employing Theorem 5.1, reduced (necessary) optimality conditions for extremals are determined by integral curves of the Hamiltonian vector field $X_h$ for the reduced Hamiltonian $h : SE(2)^2 \times (\mathfrak{se}(2)^*)^3 \to \mathbb{R}$ satisfying Hamilton's equations for $h$. Note that $\mathrm{ad}^*_{u_i} \mu_i = \begin{bmatrix} 0 & \frac{\mu_i^2 \mu_i^3}{2} & \frac{\mu_i^1 \mu_i^3}{2} \\ -\frac{\mu_i^2 \mu_i^3}{2} & 0 & \mu_i^1 \mu_i^2 \\ -\frac{\mu_i^1 \mu_i^3}{2} & -\mu_i^1 \mu_i^2 & 0 \end{bmatrix}$. Therefore by Theorem 5.1 the resulting equations are $\dot\mu_i^1 = -\mu_i^2 \mu_i^3$, $\dot\mu_i^2 = \frac{1}{2}\mu_i^1 \mu_i^3 - \sum_{j \in \mathcal{N}_i}(\Gamma_{ij})_{13}$, $\dot\mu_i^3 = -\frac{1}{2}\mu_i^1 \mu_i^2 + \sum_{j \in \mathcal{N}_i}(\widetilde\Gamma_{ij})_{23}$ together with $\dot g_i = g_i(e_1(u_i^1)^\star + e_2(u_i^2)^\star)$ for $i = 1, 2, 3$.

We show numerical simulations for the reduced necessary conditions for optimality, that is, we show the behavior of equations (10)-(11) together with the kinematic equation. We solve the corresponding initial value problem by implementing Euler's method with time step $h = 0.001$ and $N = 15000$. The agents start in an equilateral triangle with side lengths 0.5. We chose the distances $d_{ij} = 0.1$, and $\sigma_{ij} = 1$. Initial conditions are given by $g_1(0) = [\sqrt{2}/2 \; -\sqrt{2}/2 \; -1/4; \; \sqrt{2}/2 \; \sqrt{2}/2 \; 0; \; 0 \; 0 \; 1]$, $g_2(0) = [-1/\sqrt{2} \; 1/\sqrt{2} \; 1/4; \; 1/\sqrt{2} \; -1/\sqrt{2} \; 0; \; 0 \; 0 \; 1]$, $g_3(0) = [0 \; 1 \; 0; \; -1 \; 0 \; \sqrt{3}/4; \; 0 \; 0 \; 1]$, $\lambda_i = 0_{3 \times 3}$ for $i = 1, 2, 3$, $e_1 u_1^1(0) = [0 \; -5/2 \; 0; \; 5/2 \; 0 \; 0; \; 0 \; 0 \; 0]$, $e_2 u_1^2(0) = [0 \; 0 \; 5/4; \; 0 \; 0 \; 0; \; 0 \; 0 \; 0]$, $e_1 u_2^1(0) = [0 \; 2 \; 0; \; -2 \; 0 \; 0; \; 0 \; 0 \; 0]$, $e_2 u_2^2(0) = [0 \; 0 \; 2; \; 0 \; 0 \; 0; \; 0 \; 0 \; 0]$, $e_1 u_3^1(0) = [0 \; -1/2 \; 0; \; 1/2 \; 0 \; 0; \; 0 \; 0 \; 0]$, $e_2 u_3^2(0) = [0 \; 0 \; 1; \; 0 \; 0 \; 0; \; 0 \; 0 \; 0]$. In Fig. 1 we show the trajectories in the $xy$ plane (left figure) and the attitude for the agents (right figure). Fig. 2 shows the control inputs for each agent. Agents 1, 2 and 3 are given by the colors red, blue and green, respectively. Note that the blue and green agents get closer to each other in the beginning, which causes the avoidance potential to grow large (and their velocities too). The angular acceleration for the red agent is still a wave like the others agents, but with a small amplitude and oscillating along the line $u^1 = 2.5$. This is because the red agent is far from the others at their initial values, so the corresponding potential function is almost null. Observe that for the controls, there is a spike in the beginning when the agents are close, and then they relax. The control for the red agent stays close 0, which explains the low velocity and almost linear behavior of it.

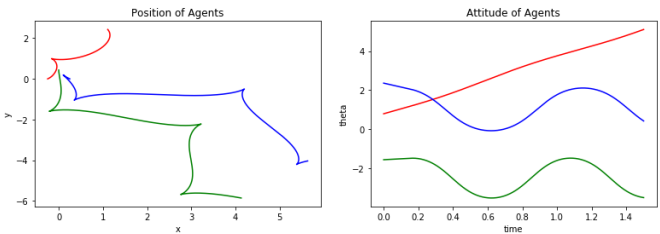

Fig. 1. Trajectories in the $xy$ (left) plane and attitude $\theta_i$ for the optimal solutions. Agents 1, 2 and 3 are given by the colors red, blue and green, respectively.

The method proposed in this work allows the construction of accurate estimators based on distance measurements, by discretizing the variational principle we proposed in Theorem 4.1 and by deriving variational integrators. This is one of the future directions of this work. Such integrators will exhibit a good performance of the energy along the motion. We will also study when such an integrator preserves the relative equilibria for the reduced system and compare the results with classical numerical methods. Reduction of sufficient conditions for optimality is also planned to be studied by using the notion of conjugate points as in [5] in a future work.

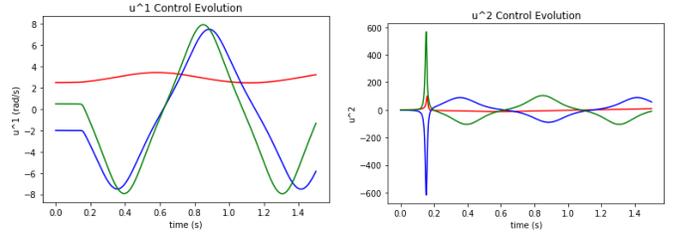

Fig. 2. Control inputs for the agents. Agents 1, 2 and 3 are given by the colors red, blue and green, respectively.